\documentclass[a4paper,12pt]{amsart}
\usepackage[letterpaper, margin=1.2in]{geometry}
\usepackage{latexsym}
\usepackage{amssymb}
\usepackage{amsmath}
\usepackage{amsfonts}
\usepackage[table]{xcolor}
\usepackage{pgfplots}
\usepackage{color}

\usepackage{txfonts,pxfonts} 
\usepackage{tikz}
\usetikzlibrary{calc,arrows}
\usetikzlibrary{calc}
\usepackage{verbatim}
\newtheorem{thm}{Theorem}[section]

\newtheorem{lem}[thm]{Lemma}

\newtheorem{cor}[thm]{Corollary}

\theoremstyle{definition}

\newtheorem{rem}{Remark}

\def\fph{\mathbb{F}_{\ph}}

\newcommand{\Z}{\mathbb Z}
\newcommand{\z}{\mathbb Z}
\newcommand{\Q}{\mathbb Q}

\newcommand{\F}{\mathbb F}
\newcommand{\fp}{\mathbb F_p}

\def\F{\mathbb{F}}

\newcommand{\p}{\mathfrak{p}}
\def\ol{\overline}
\def\al{\alpha}
\def\la{\lambda}

\def\th{\theta}
\def\md#1{\ \mbox{\rm(mod }{#1})}
\def\nph#1{N_{\ph}(#1)}
\def\npp#1{N_{\ph}^+(#1)}
\def\ph{\phi}

\newcounter{cs}
\stepcounter{cs}
\newcommand{\casos}{\begin{itemize}}
\newcommand{\fcasos}{\end{itemize}\setcounter{cs}{1}}

\newfont{\tit}{cmr12 scaled \magstep3}
\begin{document}
\title[]{On monogenity of certain  pure number fields defined by $x^{2^u\cdot 3^v}-m$}
\textcolor[rgb]{1.00,0.00,0.00}{}
\author{ Lhoussain El Fadil, and Ahmed Najim}
\address{Faculty of Sciences Dhar El Mahraz, P.O. Box  1796 Atlas-Fez , Sidi Mohamed Ben Abdellah University, Fez-- Morocco}
\email{lhoussain.elfadil@usmba.ac.ma\\najimed9@yahoo.fr}
\begin{abstract}  
 Let $K = \mathbb{Q} (\alpha) $ be a pure number field generated by   a complex root $\alpha$ of a monic irreducible polynomial  $ F(x) = x^{2^u\cdot 3^v}-m$, with $ m \neq \pm 1 $   a square free rational integer,  $u$, and $v$  two  positive  integers.   In this paper, we study the monogenity of $K$.{ The cases $u=0$ and $v=0$ have been previously studied} by the first author and Benyakkou. We prove that if $m\not\equiv 1\md4$ and $m\not\equiv \pm 1\md9$, then $K$ is monogenic. But if {$m\equiv 1\md{4}$} or $m\equiv 1\md9$ or $u=2$ and $m\equiv -1\md9$, then $K$ is not monogenic.  Some illustrating examples are given too.
\end{abstract}
\maketitle
\section{Introduction} 
Let $K=\Q(\al)$ be a number field generated by    a complex root { $\al$ } of a monic irreducible polynomial $F(x)\in
  \Z[x]$ and $\Z_K$  its ring of integers.
  It is well know that  the ring $\Z_K$ is a free $\Z$-module of rank $n=[K:\Q]$. Thus the Abelian group $\Z_K/\Z[\al]$ is finite. Its cardinal  order is called the index of $\Z[\al]$, and denoted by $(\Z_K: \Z[\al])$.
  The ring $\Z_K$ is said to have a power integral basis if it has a $\Z$-basis $(1,\th,\cdots,\th^{n-1})$ for some $\th\in \Z_K$. That is to say that $\Z_K$ is mono-generated as a ring, with a single generator $\th$. In such  a case, the field $K$ is said to
be monogenic and not monogenic otherwise. 
The problem of testing the monogenity of number fields and the construction of power integral bases has been intensively studied these last four decades, mainly by  Ga\'al, Nakahara,  Pohst, and their collaborators  (see for instance \cite{ AN,G, G19, GO, P}). 
 { In \cite{E07},   El Fadil  gave conditions for the existence of power integral bases of pure cubic fields in terms of the index form equation. In \cite{F4}, Funakura, calculated  integral bases and studied  monogenity of pure quartic fields.  In \cite{GR4},  Ga\'al and  Remete, calculated the elements of index $1$  in  pure quartic fields generated by $m^{\frac{1}{4}}$ for {$1< m <10^7$}  and $m\equiv 2,3 \md4$.  In \cite{AN6}, Ahmad, Nakahara, and Husnine proved  that  if $m\equiv 2,3 \md4$ and  $m\not\equiv \pm1\md9$, then the sextic number field generated by $m^{\frac{1}{6}}$ is monogenic.
They also showed in \cite{AN},    that if $m\equiv 1 \md4$ and $m\not\equiv \pm1\md9$, then the sextic number field generated by $m^{\frac{1}{6}}$ is not monogenic.  In \cite{E6}, based on prime ideal factorization, El Fadil showed that  if $m\equiv 1 \md4$ or $m\equiv 1\md9$, then the sextic number field generated by $m^{\frac{1}{6}}$ is not monogenic.
{Hameed and Nakahara proved that if $m\equiv 1\md4$, then the
octic number field generated by $m^{1/8}$ is not monogenic, but if $m\equiv 2,3 \md4$, 
then it is monogenic (\cite{HN8}).}   In \cite{GR17}, by applying the  explicit form of the index equation,{ Ga\'al and  Remete} obtained new deep  results on  monogenity of  number  fields generated  by $m^{\frac{1}{n}}$, with $3\le n\le 9$ and $m\neq \pm 1$ a square free integer. } In  \cite{Epr, E2r5s, E6, E6s,  E12, E24, E18,   E23k}, based on Newoton polygon's techniques, El Fadil et al. studied the monogenity of some pure  number fields. 
 The goal of this paper is to study the monogenity of pure number fields  defined by  $x^{2^u\cdot 3^v}-m$, where $m\neq \pm1$ is a square free integer, $u$, and $v$ are two natural integers. The cases { $u=0$ and  $v=0$ have been studied in \cite{Epr, CNS, BK}}. Also the case $u=1$ has been  studied by El Fadil in \cite{E23k}. 
\section{Main results}
Let $K$ be the number field generated by  a complex root{ $\al$} of   a monic  irreducible polynomial $F(x)=x^{2^u\cdot 3^v}-m$, with $m\neq \pm 1$  a square free rational integer, $u$ and $v$  two positive integers.
\begin{thm}\label{pib}
The ring $\Z[\al]$ is the ring of integers of $K$
 if and only if  $m\not\equiv 1\md4$ and  $m\not\equiv \mp 1\md9$.
\end{thm}
\begin{rem}
If $m\equiv 1\md4$ or $m\equiv 1\md9$, then $\Z[\al]$ is not integrally closed. But Theorem \ref{pib} cannot give an answer to the monogenity of $K$. The following theorems give an answer.
\end{rem}
\begin{thm}\label{npib}
If {$m\equiv 1\md{4}$} or $m\equiv 1\md9$, then $K$ is not monogenic.  
\end{thm}
{
\begin{thm}\label{npib-1}
If   $m\equiv -1\md9$ and $u=2s$ for some positive integer $s$, then $K$ is not monogenic.  
\end{thm}
\begin{cor}\label{nspe}
Let  $a\neq \pm 1$ be a square free rational integer, $u$ and $v$  two positive integers, and $k< {2^u\cdot 3^v}$ a positive integer which is coprime to $6$. Then  $F(x)=x^{2^u\cdot 3^v}-a^k$ is irreducible over $\Q$. Let $K$ be the number field defined by  a complex root { $\al$} of   a monic  irreducible polynomial $F(x)$.
\begin{enumerate}
\item
If  $a\not\equiv 1\md4$ and  $m\not\equiv \pm 1\md9$, then $K$ is monogenic.
\item
If {$a\equiv 1\md{4}$} or $a\equiv 1\md9$, then $K$ is not monogenic.  
\item 
If  $a\equiv -1\md9$ and $u=2s$ for some positive integer $s$, then $K$ is not monogenic.
\end{enumerate}
\end{cor}}
\section{Preliminaries}
{Let $K=\Q(\al)$ be a  number field generated by   a complex root $\al$  of a monic irreducible polynomial $F(x)\in \Z[x]$, $\Z_K$ its ring of integers, and  $ind(\al)=(\Z_K: \Z[\al])$ the index of $\Z[\al]$ in $\Z_K$.
  For a rational prime integer $p$, if $p$ does not divide $(\mathbb{Z}_K : \mathbb{Z}[\alpha])$, then a well known theorem of Dedekind says that the  factorization of  $ p \mathbb{Z}_K$ can be derived directly  from the factorization of $\overline{F(x)}$ in $\F_p[x]$. Namely,  $p\Z_K=\prod_{i=0}^r \p_i^{l_i}$, where every $\p_i=p\Z_K+\phi_i(\al)\Z_K$ and $\ol{F(x)}=\prod_{i=1}^r \ol{\ph_i(x)}^{l_i}$ modulo $p$ is the factorization of $\ol{F(x)}$ into powers of monic irreducible coprime polynomials of $\F_p[x]$. So, $f(\p_i)={\mbox{deg}}(\phi_i)$ is the residue degree of $\p_i$ (see \cite[ Chapter I, Proposition 8.3]{Neu}).
  In order to apply this  theorem in an effective way, one needs a
criterion to test  whether  $p$ divides  or not the index $(\Z_K:\Z[\al])$.  For a number field $K$ generated by  a complex root $\al$ of a monic irreducible  polynomial $F(x)\in \Z[x]$ and a rational prime integer $p$, let $\overline{F(x)}=\prod_{i=1}^r\overline{\ph_i(x)}^{l_i}\md{p}$  be the factorization of   $\overline{F(x)}$ in $\F_p[x]$, where the polynomials $\ph_i\in\Z[x]$ are monic with  their reductions are irreducible over $\F_p$ and GCD$(\overline{\ph_i},\overline{\ph_j})=1$ for every $i\neq j$. If we set
				$M(x)=\cfrac{F(x)-\prod_{i=1}^r{\ph_i}(x)^{l_i}}{p}$. In $1878$,  Dedekind proved that $M(x)\in \Z[x]$ and  the  well known Dedekind's criterion:
				\begin{thm}\label{Ded}$($\cite[Theorem 6.1.4]{Co} and \cite{R}$)$\\
			The following statements are equivalent:
				\begin{enumerate}
					\item[1.]
					$p$ does not divide the index $(\Z_K:\Z[\al])$.
					\item[2.]
					For every $i=1,\dots,r$, either $l_i=1$ or $l_i\ge 2$ and $\overline{\ph_i(x)}$ does not divide $\overline{M(x)}$ in $\F_p[x]$.
				\end{enumerate}
		\end{thm}
When Dedekind's criterion fails, that is,  $p$ divides the index $(\z_K:\z[\th])$  for every primitive element $\th\in \Z_K$ of $K$, {  then} it is not possible to obtain the prime ideal factorization of $p\z_K$ by applying { Dedekind's theorem}.
In 1928, Ore developed   an alternative approach
for obtaining the index $(\z_K:\z[\alpha])$, the
absolute discriminant, and the prime ideal factorization of the rational primes in
a number field $K$ by using Newton polygons (see \cite{MN, O}).
For the convenience of the reader, as it is necessary for the proof of our main results,  the content of this section is copied from \cite{E12}.\\
  We start by recalling  some fundamental facts about Newton polygons  applied in algebraic number theory. For more details, we refer to \cite{El, EMN, GMN}. 
 For a  prime integer  $p$ and for a  monic polynomial 
$\phi\in \z[x]$ {whose reduction} is irreducible  in
$\fp[x]$, 
let $\fph$ be 
the field $\frac{\fp[x]}{(\overline{\phi})}$. For any
monic polynomial  $F(x)\in \z[x]$, upon  the Euclidean division
 by successive powers of $\ph$, we  expand $F(x)$ as
$F(x)=\sum_{i=0}^la_i(x)\phi(x)^{i}$, called    the $\phi$-expansion of $F(x)$
 (for every $i$, $deg(a_i(x))<
deg(\phi)$). 
The $\ph$-Newton polygon of $F(x)$ with respect to $p$, is the lower boundary of the convex envelope of the set of points $\{(i,\nu_p(a_i(x))),\, a_i(x)\neq 0\}$ in the Euclidean plane, which is denoted by $\nph{F}$. Let $S_1,\, S_2,\dots, S_t$ be the sides of $\nph{F}$. For every side $S$ of $\nph{F}$, the length of $S$, denoted by $l(S)$, is the length of its projection to the $x$-axis, its height, denoted by $H(S)$, is the length of its projection to the $y$-axis. Let $\la=H(S)/l(S)$, then $-\la$ is the slope of $S$. If  $\la\neq 0$, then  $\la=h/e$ with $e$ and $h$  two positive coprime integer.  Notice that $e=l(S)/d$, called the ramification index of $S$ and $h=H(S)/d$, where   $d=$gcd$(l(S),H(S))$ is called the degree of $S$. Thus $\nph{F}$ is the join of its different sides ordered by increasing slopes, which we can express by $\nph{F}=S_1+ S_2+\dots+ S_t$.
 The principal $\ph$-Newton polygon of $F(x)$ ,
 denoted by $\npp{F}$, is the part of the  polygon $\nph{F}$, which is  determined by  all sides of negative  slopes of $\nph{F}$.
  For every side $S$ of {$\npp{F}$}, with initial point $(s, u_s)$ and length $l$, and for every 
$i=0, \dots,l$, we attach   the following
 residue coefficient $c_i\in\fph$ as follows:
 $$c_{i}=
\left
\{\begin{array}{ll} 0,& \mbox{ if } (s+i,{\it u_{s+i}}) \mbox{ lies strictly
above } S
 ,\\
\left(\dfrac{a_{s+i}(x)}{p^{{\it u_{s+i}}}}\right)
\,\,
\md{(p,\phi(x))},&\mbox{ if }(s+i,{\it u_{s+i}}) \mbox{ lies on }S.
\end{array}
\right.$$
where $(p,\phi(x))$ is the maximal ideal of $\z[x]$ generated by $p$ and $\ph$.\\
Let  $\la=-h/e$ be the slope of $S$, where  $h$ and $e$ are two positive coprime integers. Then $d=l/e$ is the degree of $S$.  Notice that, the points  with integer coordinates lying on $S$ are exactly $(s,u_s),(s+e,u_{s}-h),\cdots, (s+de,u_{s}-dh)$. Thus, if $i$ is not a multiple of $e$, then 
$(s+i, u_{s+i})$ does not lie on $S$, and so $c_i=0$. Let
$F_S(y)=t_dy^d+t_{d-1}y^{d-1}+\cdots+t_{1}y+t_{0}\in\fph[y]$, called  
the residual polynomial of $F(x)$ associated to the side $S$, where for every $i=0,\dots,d$,
 $t_i=c_{ie}$.
\begin{rem}
\begin{enumerate}
\item
Notice that, since  $(s, u_{s})$ and $(s+l, u_{s+l})$ lie on $S$, we conclude that $t_dt_0\neq 0$ in $\fph$,  and so deg$(F_S)=d$ and $F_S(0)\neq 0$.
\item
 Notice also that if $\nu(a_{s}(x))=0$, $\la=0$, and $\ph=x$, then $\fph=\F_p$ and for every $i=0,\dots,l$, $c_i=\overline{{a_{s+i}}} \md{p}$. Thus this notion of residual coefficient generalizes the reduction modulo the maximal ideal $(p)$ and $F_S(y)\in\F_p[y]$ coincides with the reduction of $F(x)$ modulo the maximal ideal $(p)$.
\end{enumerate}
    \end{rem}
    Let $\npp{F}=S_1+\dots + S_t$ be the { principal} $\ph$-Newton polygon of $F$ with respect to $p$. We say that $F(x)$ is a $\ph$-regular polynomial with respect to $p$,
    if for every  $i=1,\dots,t$, $F_{S_i}(y)$ is square free in $\fph[y]$. We say that $F$ is a $p$-regular polynomial if $F$ is  a $\ph_i$-regular polynomial with respect to $p$ for every $i=1,\dots,r$, for some monic polynomials $\ph_1,\dots,\ph_r$ in $\Z[x]$, with $\overline{\ph_1},\dots,\overline{\ph_r}$ are pairwise coprime irreducible polynomials and  $\overline{F(x)}=\prod_{i=1}^r\overline{\ph_i}^{l_i}$ is the factorization of $\overline{F(x)}$ in $\F_p[x]$.\\\\
The  theorem of Ore plays a key role for proving our main theorems:\\
  Let $\ph\in\Z[x]$ be a monic polynomial, with $\overline{\ph(x)}$  irreducible in $\F_p[x]$. As defined in \cite[Def. 1.3]{EMN},   the $\ph$-index of $F(x)$, denoted  $ind_{\ph}(F)$, is  deg$(\ph)$ multiplied by  the number of points with natural integer coordinates that lie below or on the polygon $\npp{F}$, strictly above the horizontal axis, {and strictly beyond the vertical axis} (see Figure 1).
 \begin{figure}[htbp] 
\begin{tikzpicture}[x=1cm,y=0.5cm]
\draw[latex-latex] (0,6) -- (0,0) -- (10,0) ;
\draw[thick] (0,0) -- (-0.5,0);
\draw[thick] (0,0) -- (0,-0.5);
\node at (0,0) [below left,blue]{\footnotesize $0$};
\draw[thick] plot coordinates{(0,5) (1,3) (5,1) (9,0)};
\draw[thick, only marks, mark=x] plot coordinates{(1,1) (1,2) (1,3) (2,1)(2,2)     (3,1)  (3,2)  (4,1)(5,1)  };
\node at (0.5,4.2) [above  ,blue]{\footnotesize $S_{1}$};
\node at (3,2.2) [above   ,blue]{\footnotesize $S_{2}$};
\node at (7,0.5) [above   ,blue]{\footnotesize $S_{3}$};
\end{tikzpicture}
\caption{    \large  $\npp{F}$.}
\end{figure}
 Now assume that $\overline{F(x)}=\prod_{i=1}^r\overline{\ph_i}^{l_i}$ is the factorization of $\overline{F(x)}$ in $\F_p[x]$,  {where  $\ph_1,\dots, \ph_r$ are monic polynomials lying in $\Z[x]$ and $\ol{\ph_1},\dots, \ol{\ph_r}$ are pairwise coprime irreducible polynomials over $\F_p$}.
For every $i=1,\dots,r$, let  $N_{\ph_i}^+(F)=S_{i1}+\dots+S_{ir_i}$ be the principal part of the $\ph_i$-Newton polygon of $F$ with respect to $p$. For every {$j=1,\dots, r_i$},  let $F_{S_{ij}}(y)=\prod_{s=1}^{s_{ij}}\psi_{ijs}^{a_{ijs}}(y)$ be the factorization of $F_{S_{ij}}(y)$ { into powers of monic irreducible polynomials } of $\F_{\ph_i}[y]$. 
  Then we have the following  theorem of Ore (see \cite[Theorem 1.7 and Theorem 1.9]{EMN}, \cite[Theorem 3.9]{El}, and \cite{MN}):
 \begin{thm}\label{ore} (Theorem of Ore)
 \begin{enumerate}
 \item
  {$$\nu_p((\z_K:\z[\al]))\ge \sum_{i=1}^r ind_{\ph_i}(F).$$} 
  The equality holds if $F(x)$ is $p$-regular.
\item
If  $F(x)$ is $p$-regular, then
$$p\Z_K=\prod_{i=1}^r\prod_{j=1}^{r_i}
\prod_{s=1}^{s_{ij}}\p^{e_{ij}}_{ijs},$$ where $e_{ij}$ is the ramification index
 of the side $S_{ij}$ and $f_{ijs}=\mbox{deg}(\ph_i)\times \mbox{deg}(\psi_{ijs})$ is the residue degree of $\p_{ijs}$ over $p$ for every $i=1,\dots,r$, $j=1,\dots,r_i$, and $s=1,\dots, s_{ij}$.
 \end{enumerate}
\end{thm}
\begin{cor}\label{indore}
{ Under the assumptions above Theorem \ref{ore}, if for every $i=1,\dots,r$, $l_i=1$ or $N_{\ph_i}(F)=S_i$ has a single side of height $1$, then $\nu_p((\z_K:\z[\al]))=0$}.
\end{cor}
The following lemma allows to evaluate  the $p$-adic valuation of the binomial coefficient $\binom{p^r}{j} $. Its proof will appear in \cite{E2r5s}.
\begin{lem} \label{binomial}
	Let $p$ be a rational prime integer and $r$ be a positive integer. Then $ \nu_p\left(\binom{p^r}{j}\right)  =  r - \nu_p(j)$
	for any integer $j= 1,\dots,p^r-1 $. 
\end{lem}
In \cite{GMN}, Guàrdia, Montes, and  Nart introduced  the notion of $\phi$-admissible expansion used in order to treat some special cases when the $\phi$-expansion is hard to calculate. Let
\begin{equation}
\label{eq1}
F(x)=\sum_{i=0}^nA_i'(x)\phi(x)^i,\quad A_i'(x)\in \mathbb{Z}[x],
\end{equation}
be a $\phi$-expansion of $F(x)$, not necessarily the $\phi$-expansion (deg$(A_i')$ is not necessarily less than deg$(\ph)$). Take $u_i'=\nu_p(A_i'(x))$, for all $i=0,\dots, n$, and let $N'$ be the lower boundary of the convex envelope of the set of points $\{(i,u_i')\,\mid\,0\leq i\leq n,\,u_i'\neq\infty\}$ and $N'^+$ its principal part. To any $i=0,\dots,n$, we attach the residue coefficient as follows:
$$
c_i'=\left\{
\begin{array}{ll}
0,&\text{if }(i,u_i')\text{ lies above }N',\\
\left(\frac{A_i'(x)}{p^{u_i'}}\right)(\mod(p,\phi(x))),& \text{if }(i,u_i')\text{ lies on }N'.
\end{array}
\right.
$$
Likewise, for any side $S$ of $N'^+$, we can define the residual polynomial attached to $S$ and denoted $R_{\lambda}'(F)(y)\,($similar to the residual polynomial $R_{\lambda}(F)(y)$ from the $\phi$-adic expansion$)$. We say that the $\phi$-expansion $(\ref{eq1})$ is admissible if $c_i'\neq 0$ for each abscissa $i$ of a vertex of $N'$. For more details, we refer to \cite{GMN}.
\begin{lem}{$($\cite[Lemma 1.12]{GMN}$)$}\\
	If a $\phi$-expansion of $F(x)$ is admissible, then $N'^+=\npp{F}$ and $c_i'=c_i$. In particular, for any side $S$ of $N'^+$ we have $R_{\lambda}'(F)(y)=R_{\lambda}(F)(y)$ up to {multiplication} by a nonzero coefficient of $\fph$.
\end{lem}
The following lemma  allows to determine the $\ph$-Newton polygon of $F(x)$. Its proof will appear in \cite{E3u7v}.
\begin{lem}\label{NP}
Let $F(x)=x^n-m\in \Z[x]$ be an irreducible polynomial and $p$ a prime integer which divides $n$ and does not divide $m$. Let $n=p^rt$ in $\Z$ with $p$ does not divide $t$. Then $\ol{F(x)}=\ol{(x^t-m)}^{p^r}$. Let $v=\nu_p(m^p-m)$ and $\ph\in \Z[x]$ be a monic polynomial, whose reduction modulo $p$ divides $\ol{F(x)}$. 
 \begin{enumerate}
 \item
 If $\nu_p(m^{p-1}-1)\le r$, then $\npp{F}$ is the lower boundary of the  convex envelope of the set of the points $\{(0,v)\}\cup \{(p^j,r-j), \, j=0,\dots,r\}$.
 \item
 If $\nu_p(m^{p-1}-1)\ge r+1$, then $\npp{F}$ is the lower boundary of the  convex envelope of the set of the {points} $\{(0,V)\}\cup \{(p^j,r-j), \, j=0,\dots,r\}$ for some integer $V\ge r+1$.
 \end{enumerate}
\end{lem}}
\section{Proofs of main results}
\begin{proof} of Theorem \ref{pib}.\\
The proof of Theorem \ref{pib} can be done by using Dedekind's criterion as it was shown in the proof of \cite[Theorem $6.1$]{HS}. But as the other results are based on Newton polygon's techniques, let us use theorem of index with "if and only if" as it is given in  \cite[Theorem 4.18]{GMN}, which  says that: $\nu_p(\Z_K:\Z[\al])=0$ if and only if $ind_1(F)=0$, where $ind_1(F)$ is the index given in Theorem \ref{ore}.  Since $\triangle(F)  =\mp (2^{u}3^{v})^{2^{u}3^{v}} m ^{2^{u}3^{v}}$, then by the  formula  $\nu_p(\triangle(F))=2\nu_p(ind(F))+\nu_p(d_K)$, where $d_K$ is the absolute discriminant of $K$ and $ ind(F)=(\Z_K: \Z[\al])$, we conclude that $\mathbb{Z}[\alpha]$ is integrally closed if and only if   $p$ does not divide $(\Z_K: \Z[\al])$ for every rational prime integer $p$ dividing $6m$.  Let $p$ be a rational prime dividing $m$, then $ F(x) \equiv \ph^{2^{u}\cdot 3^{v}}  (\ mod \ p )$, where $\phi = x$. As $m$ is a square free integer, the $\phi$-principal Newton polygon with respect to $\nu_p$,  $ \npp {F} = S $ has a single side of height $\nu_p(m)$. As $l(S)=2^{u}\cdot 3^{v}$, $ind_\ph(F)=0$ if and only if the height of $S$ equals $1$, which  means $\nu_p(m)=1$. 
  It follows that the unique prime candidates to divide the index $(\mathbb{Z}_K:\mathbb{Z}[\alpha])$ are $2$ and $3$.\\
 For $p=2$  and $2$ does not divide $m$,  let $\ph\in \Z[x]$ be a monic polynomial, whose reduction is an irreducible factor of  $(x^{3^v}-1)$ in $\F_2[x]$.
 As $l(\npp{F})=2^{u}\ge 2$, $ind_\ph(F)=0$ if and only if  $\npp{F}$ has a single side of height  $1$, which  means by Lemma \ref{NP} that $ \nu_2(1- m)=1=1$; $m\equiv 3\md4$.\\  
 Similarly, for  $p=3$  and $3$ does not divide $m$,  let  $\ph\in \Z[x]$ be a monic polynomial, whose reduction is an irreducible factor of  $(x^{2^u}-1)$ in $\F_3[x]$.
 Again  as $l(\npp{F})=3^{v}\ge 2$, $ind_\ph(F)=0$ if and only if  $\npp{F}$ has a single side of height  $1$, which  means by Lemma \ref{NP} that $ \nu_3(m^2- 1)=1=1$; $m\not\equiv \mp1\md9$.  
  \end{proof} 
 \smallskip
	
	{The index of a field $K$ is defined by $i(K)=gcd\{(\mathbb{Z}_K:\mathbb{Z}[\theta])\mid K=\mathbb{Q}(\theta) \mbox{ and } \theta\in \Z_K \}$. A rational prime $p$ dividing $i(K)$ is called a prime common index divisor of $K$. If $\mathbb{Z}_K$ has a power integral basis, then $i(K)=1$. Therefore a field having a prime common index divisor is not monogenic. The existence of prime common index divisors was first established in $1871$ by Dedekind who exhibited examples in fields of third and fourth degrees, for example, he considered the cubic field $K$ defined by $F(x)=x^3-x^2-2x-8$ and he showed that the prime $2$ splits completely. So, if we suppose that $K$ is monogenic, then we would be able to find a cubic polynomial generating $K$, that splits completely into distinct polynomials of degree $1$ in $\mathbb{F}_2[x]$. Since there are only $2$ distinct polynomials of degree $1$ in $\mathbb{F}_2[x]$, this is impossible. Based on these ideas and using Kronecker's theory of algebraic numbers, Hensel  gave a necessary and sufficient condition on the so-called "index divisors" for any prime integer  $p$ to be  a prime common index divisor \cite{He2}. $($For more details see \cite{HS}$)$}.
	{For the proof of Theorem $\ref{npib}$, we need the following lemma, which  characterizes the prime common index divisors of $K$.  We need to use only  one way, which is an immediate consequence of  Dedekind's theorem}.
	 	\begin{lem}\label{index}$($\cite[Theorem 2.2]{HS}$)$\\
		Let $p$ be a rational prime integer and $K$ be a number field. For every positive integer $f$, let $\mathcal{P}_f$ be the number of distinct prime ideals of $\mathbb{Z}_K$ lying above $p$ with residue degree $f$ and $\mathcal{N}_f$ the number of monic irreducible polynomials of $\mathbb{F}_p[x]$ of degree $f$.  Then $p$ is a prime common index divisor of $K$ if and only if $\mathcal{P}_f>\mathcal{N}_f$ for some positive integer $f$.
	\end{lem}
	\begin{rem}
		{As it was shown in the proof of Theorem \ref{pib} that: the unique prime candidates to be a prime common index divisors of $K$ are $2$ and $3$, because if $p\not\in\{2,3\}$, then $p$ does not divide the index $(\mathbb{Z}_K:\mathbb{Z}[\alpha])$, and so the factorization of $p\mathbb{Z}_K$ is analogous to the factorization of $x^{2^u\cdot 3^v}-m$ in $\mathbb{F}_p[x]$.}
	\end{rem}
	\begin{rem}\label{remore}
		{In order to prove  Theorem \ref{npib}, we don't need to determine the factorization of $p\Z_K$ explicitly. But according to Lemma \ref{index}, we need only to show that $\mathcal{P}_f>\mathcal{N}_f$ for an adequate positive integer  $f$. So in practice the second point of Theorem \ref{ore}, could be replaced by the following:
			If  $l_i=1$ or $d_{ij}=1$ or $a_{ijk}=1$ for some $(i,j,k)$ according to notation of    Theroem \ref{ore}, then $\psi_{ijk}$ provides  a prime ideal $\p_{ijk}$ of $\Z_K$ lying above $p$ with residue degree  {$f_{ijk}=m_i\times t_{ijk}$}, where  $t_{ijk}=$deg$(\psi_{ijk})$ and $p\Z_K=\p_{ijk}^{e_ij}I$, where the factorization of the ideal $I$ can be derived from the other factors of each residual polynomials of $F(x)$.}
	\end{rem}
 \begin{proof} of Theorem \ref{npib}.
 \begin{enumerate}
 \item
  If    $m\equiv 1\md2$, then $\overline{F(x)}=\overline{(x^{3^v}-1)^{2^u}}=\overline{((x^{3}-1)U(x))}^{2^u}=\overline{(x^2+x+1)U(x)}^{2^u}$  in $\F_2[x]$ for a monic polynomial $U(x)\in\Z[x]$.
  
{  Let $\ph_2=x-1$, $\ph_2=x^2+x+1$, and $v_2=\nu_2(1-m)$.} 
  \begin{enumerate}
 \item
{If $v_2=2$, then $N_{\ph_2}^+(F)=S$ has a single side of degree $d=2$. By using $F(x) = (\ph_1U(x))^{2^u}+\cdots+1-m$, we have 
 $F_{S}(y)=t^2y^2+ty+1$, where $t\equiv \ph_1U(x)\md{2,\ph_2}$ si a nonzero element of $\F_{\ph_2}$. Since $x^{3^k}-1$ is separable over $\F_2$, $\ol{\ph_2}$ does not divide $\ol{\ph_1U(x)}$ in $\F_2[x]$, and so  $t$ is a nonzero element of  $\F_{\ph_2}$. It follows $F_{S}(y)=(ty-x)(ty-x^2)$ in $\F_{\ph_2}[y]$. Thus $\ph_2$ provides $2$ distinct prime ideals of $\Z_K$ lying above $2$ with residue degree $2$ each.}
 If  $u\ge 2$ and  $v_2=3$, then by{ Lemma \ref{NP}}, $N_{\ph_i}^+(F)$ has  two sides $S_{i1}$ and $S_{i2}$ joining the point $(0,3)$, $(2^{u-1},1)$, and  $(2^u,0)$ (see $FIGURE\, 2$).    Thus $S_{i1}$ is a side of  degree $2$  and $S_{i2}$ is a side of  degree $1$ for every $i=1,2$. Since  $F_{S_{11}}(y)=y^2+y+1$ is irreducible over $\F_{\ph_1}\simeq \F_2$ and  $F_{S_{22}}(y)$ is of degree $1$. By Remark \ref{remore}, every  $\ph_i$ provides at least prime ideal $\p_{ii}$ of $\Z_K$ associated to the side $S_{ii}$, with  $f_{11}=$deg$(\ph_1)\cdot$deg$(F_{S_{11}}=1\cdot 2=2$ and $f_{22}=$deg$(\ph_2)\cdot$deg$(F_{S_{22}})=2\cdot 1=2$ are  the residue degrees of $\p_{11}$ and $\p_{22}$. Thus there are at least $2$ distinct prime ideals of $\Z_K$ lying above $2$ with residue degree $2$ each. As $x^2+x+1$ is the unique monic irreducible polynomial of degree $2$ in $\F_2[x]$,  by Lemma \ref{index}, $2$ divides $i(K)$ and $K$ is not monogenic. 
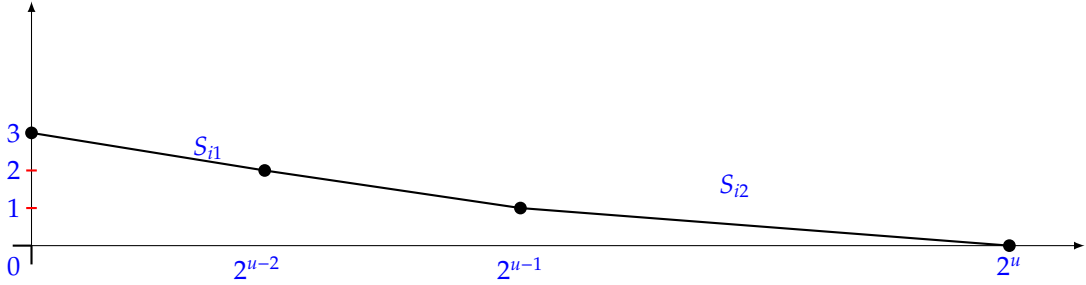
\begin{figure}[htbp] 
\begin{tikzpicture}[x=0.5cm,y=0.5cm]
\draw[latex-latex] (0,6.5) -- (0,0) -- (28,0) ;
\draw[thick] (0,0) -- (-0.5,0);
\draw[thick] (0,0) -- (0,-0.5); 
\draw[thick,red] (-2pt,1) -- (2pt,1);
\draw[thick,red] (-2pt,2) -- (2pt,2);
\draw[thick,red] (-2pt,3) -- (2pt,3);
\node at (0,0) [below left,blue]{\footnotesize  $0$};
\node at (6,0) [below ,blue]{\footnotesize  $2^{u-2}$};
\node at (13,0) [below ,blue]{\footnotesize  $2^{u-1}$};
\node at (26,0) [below ,blue]{\footnotesize  $2^{u}$};
\node at (0,1) [left ,blue]{\footnotesize  $1$};
\node at (0,2) [left ,blue]{\footnotesize  $2$};
\node at (0,3) [left ,blue]{\footnotesize  $3$};
\draw[thick,mark=*] plot coordinates{(0,3)(6.2,2)(13,1) (26,0)};
\node at (4,2) [above right  ,blue]{\footnotesize  $S_{i1}$};
\node at (18,1) [above right  ,blue]{\footnotesize  $S_{i2}$};
\end{tikzpicture}
\caption{    \large   $N_{\ph_i}^+(F)$ for  $v_2=3$.\hspace{5cm}}
\end{figure}
 \item
 If  $u\ge 2$ and  $v_2\ge 4$, then by Lemma \ref{NP}, $N_{\ph_i}^+(F)$ has at least $3$ sides for which the last two sides $S_{i1}$ and $S_{i2}$ are of height $1$ for every $i=1,2$ (see $FIGURE 3$).    Thus,  $F_{S_{21}}(y)$ and $F_{S_{22}}(y)$ are of degree $1$. By Remark \ref{remore},   $\ph_2$ provides at least two  prime ideals $\p_{21}$ and $\p_{22}$ of $\Z_K$ lying above $2$ with  residue degree $f_{2i}=$deg$(\ph_2)\cdot$deg$(F_{S_{2i}}=2$ each. As $x^2+x+1$ is the unique monic irreducible polynomial of degree $2$ in $\F_2[x]$,  by Lemma \ref{index}, $2$ divides $i(K)$ and $K$ is not monogenic. 
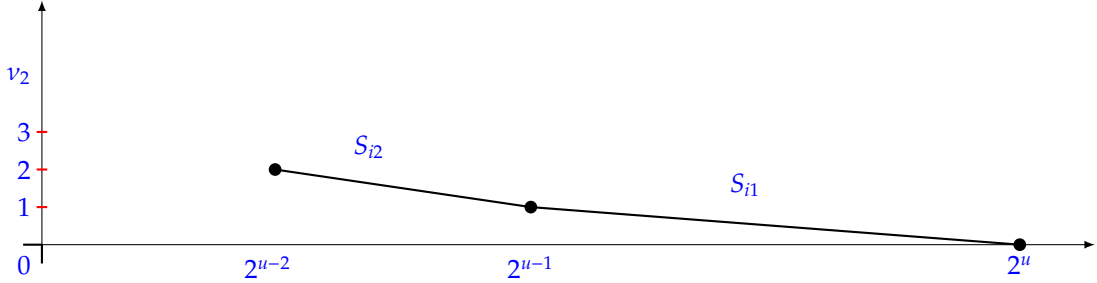
\begin{figure}[htbp] 
\begin{tikzpicture}[x=0.5cm,y=0.5cm]
\draw[latex-latex] (0,6.5) -- (0,0) -- (28,0) ;
\draw[thick] (0,0) -- (-0.5,0);
\draw[thick] (0,0) -- (0,-0.5); 
\draw[thick,red] (-2pt,1) -- (2pt,1);
\draw[thick,red] (-2pt,2) -- (2pt,2);
\draw[thick,red] (-2pt,3) -- (2pt,3);
\node at (0,0) [below left,blue]{\footnotesize  $0$};
\node at (6,0) [below ,blue]{\footnotesize  $2^{u-2}$};
\node at (13,0) [below ,blue]{\footnotesize  $2^{u-1}$};
\node at (26,0) [below ,blue]{\footnotesize  $2^{u}$};
\node at (0,1) [left ,blue]{\footnotesize  $1$};
\node at (0,2) [left ,blue]{\footnotesize  $2$};
\node at (0,3) [left ,blue]{\footnotesize  $3$};
\node at (0,4.5) [left ,blue]{\footnotesize  $\nu_2$};
\draw[thick,mark=*] plot coordinates{(6.2,2)(13,1) (26,0)};
\node at (8,2) [above right  ,blue]{\footnotesize  $S_{i2}$};
\node at (18,1) [above right  ,blue]{\footnotesize  $S_{i1}$};
\end{tikzpicture}
\caption{    \large   $N_{\ph_i}^+(F)$ for  $v_2\ge 4$.\hspace{5cm}}
\end{figure}
 \end{enumerate}
\item
If $m\equiv 1\md9$, then $\overline{F(x)}=(x^{2^u}-1)^{3^v}=((x-1)(x+1)U(x))^{3^v}$ in $\F_3[x]$ for a monic polynomial  $U(x)\in \F_3[x]$. Let $\ph_1=x-1$, $\ph_2=x+1$, and $v_3=\nu_3(1-m)$.  If  $v_3\ge 2$, then by Lemma \ref{NP}, $N_{\ph_i}^+(F)$ has at least $2$ sides of which the last two sides $S_{i1}$ and $S_{i2}$ are of height $1$ each for every $i=1,2$ (see $FIGURE 4$ and $FIGURE 5$).
Thus  $F_{S_{ij}}(y)$ is of degree $1$ for every $i,j=1,2$. By Remark \ref{remore}, every  $\ph_i$ provides at least $2$  prime ideals $\p_{ij}$ of $\Z_K$ lying above $3$ with  residue degree $f_{ij}=$deg$(\ph_i)\cdot$deg$(F_{S_{ij}})=1$ for every $i,j=1,2$, and so there are at least $4$ prime ideals of $\Z_K$ lying above $3$ with  residue degree $1$ each. As  there is only $3$ monic irreducible polynomial of degree $1$ in $\F_3[x]$,  by Lemma \ref{index}, $3$ divides $i(K)$ and $K$ is not monogenic. 
\begin{figure}[htbp] 
\begin{tikzpicture}[x=0.5cm,y=0.5cm]
\draw[latex-latex] (0,6.5) -- (0,0) -- (28,0) ;
\draw[thick] (0,0) -- (-0.5,0);
\draw[thick] (0,0) -- (0,-0.5); 
\draw[thick,red] (1,-2pt) -- (1,2pt);
\draw[thick,red] (3,-2pt) -- (3,2pt);
\draw[thick,red] (9,-2pt) -- (9,2pt);
\draw[thick,red] (-2pt,1) -- (2pt,1);
\draw[thick,red] (-2pt,2) -- (2pt,2);
\draw[thick,red] (-2pt,3) -- (2pt,3);
\node at (0,0) [below left,blue]{\footnotesize  $0$};
\node at (3,0) [below ,blue]{\footnotesize $3^{v-2}$};
\node at (9,0) [below ,blue]{\footnotesize  $3^{v-1}$};
\node at (27,0) [below ,blue]{\footnotesize  $3^{v}$};
\node at (0,1) [left ,blue]{\footnotesize  $1$};
\node at (0,2) [left ,blue]{\footnotesize  $2$};
\node at (0,3) [left ,blue]{\footnotesize  $v_3$};
\draw[thick,mark=*] plot coordinates{(0,3) (3,2) (9,1) (27,0)};
\node at (6,1.3) [above right  ,blue]{\footnotesize  $S_{i2}$};
\node at (16,0.5) [above right  ,blue]{\footnotesize  $S_{i1}$};
\end{tikzpicture}
\caption{    \large   $N_{\ph_i}^+({F})$ for  $v\ge 2$ and $v_3=3$.\hspace{5cm}}
\end{figure}
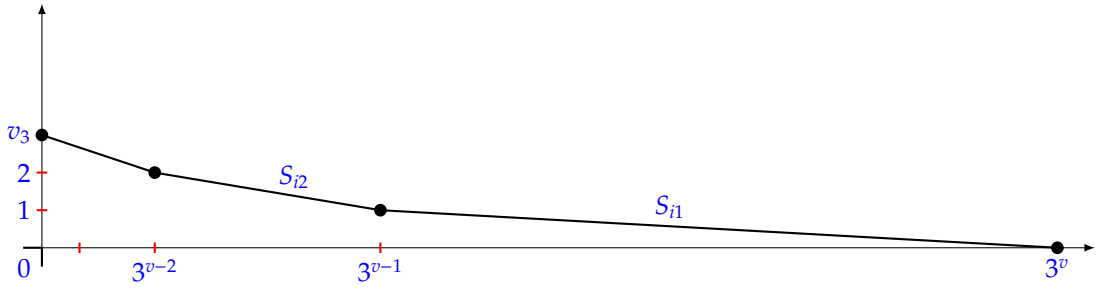
\newpage 
\begin{figure}[htbp] 
\begin{tikzpicture}[x=0.5cm,y=0.5cm]
\draw[latex-latex] (0,6.5) -- (0,0) -- (28,0) ;
\draw[thick] (0,0) -- (-0.5,0);
\draw[thick] (0,0) -- (0,-0.5); 
\draw[thick,red] (1,-2pt) -- (1,2pt);
\draw[thick,red] (3,-2pt) -- (3,2pt);
\draw[thick,red] (9,-2pt) -- (9,2pt);
\draw[thick,red] (-2pt,1) -- (2pt,1);
\draw[thick,red] (-2pt,2) -- (2pt,2);
\draw[thick,red] (-2pt,3) -- (2pt,3);
\node at (0,0) [below left,blue]{\footnotesize  $0$};
\node at (3,0) [below ,blue]{\footnotesize $1$};
\node at (9,0) [below ,blue]{\footnotesize  $3$};
\node at (0,1) [left ,blue]{\footnotesize  $1$};
\node at (0,2) [left ,blue]{\footnotesize  $2$};
\node at (0,4) [left ,blue]{\footnotesize  $v_3$};
\draw[thick,mark=*] plot coordinates{(0,4) (3,1) (9,0)};
\node at (2,2) [above right  ,blue]{\footnotesize  $S_{i2}$};
\node at (5.5,0.8) [above right  ,blue]{\footnotesize  $S_{i1}$};
\end{tikzpicture}
\caption{    \large   $N_{\ph_i}^+({F})$ for  $v=1$ and $v_3\ge 3$.\hspace{5cm}}
\end{figure}
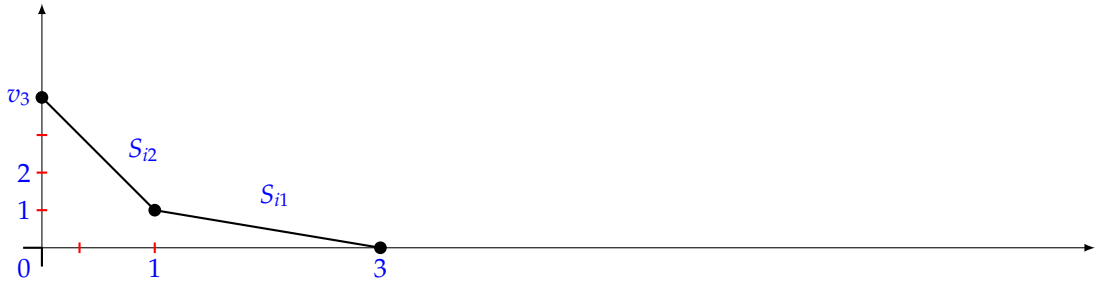
\end{enumerate}
\end{proof}
\begin{proof} of Theorem \ref{npib-1}.
 Assume that  $m\equiv -1\md9$ and $u=2s$ for some positive integer $s$, then $\overline{F(x)}=\overline{((x^{4}+1)U(x))}^{{3^v}}=((x^{2}+x-1)(x^2-x-1))^{3^v}$ in $\F_3[x]$ for monic polynomial $U\in \Z[x]$ such that ${\overline{x^{4}+1}}$ and ${\overline{U(x)}}$ are coprime in $\F_3[x]$. Let $\ph_1=x^2+x-1$, $\ph_2=x^2-x-1$, and $v_3=\nu_3(1-m)$. 
 By Lemma \ref{NP}, if $v_2\ge 2$, then $N_{\ph_i}^+(F)$ has at least two sides $S_{i1}$ and $S_{i2}$ with height $1$ each. Thus $F_{S_{ij}}(y)$ is irreducible over $\F_{\ph_i}$ as it is of degree $1$ for every $i,j=1,2$. By Remark \ref{remore}, every factor $\ph_i$ provides at least two  distinct prime ideals of $\Z_K$ lying above $3$ with residue degree $f=2$ each. Thus there are at least four  distinct prime ideals of $\Z_K$ lying above $3$ with residue degree $f=2$ each. As  $x^2+1$, $x^2+x-1$, and  $x^2-x-1$ are  the unique monic irreducible polynomials of degree $2$ in $\F_3[x]$,  by Lemma \ref{index}, $3$ divides $i(K)$, and so $K$ is not monogenic.
\end{proof}
\begin{proof} of Corollary \ref{nspe}.\\
Since gcd$(k,6)=1$, let $(x,y)\in \z^2$ be the unique solution of the equation $k\cdot x- 2^u\cdot 3^v\cdot y=1$ and $\th=\frac{\al^x}{a^y}$. Then $\th^{2^u\cdot 3^v}=a$, and so $g(x)=x^{2^u\cdot 3^v}-a$ is the minimal polynomial of $\th$ over $\Q$, $\th\in \z_K$ is a primitive element of $K$. Since $a$ is square free, we can apply Theorems \ref{pib} and \ref{npib}.
\end{proof}
\section{Examples}
 Let $K=\Q(\al)$ be the pure  number fields generated by $\al$ a root of a monic irreducible polynomial  $F(x)=x-^{2^u\cdot 3^v}-m$ with $u$ and $v$ are two positive integers.
\begin{enumerate}
    \item 
    Let  $F(x)=x^{36}-11664\in \Z[x]$. Since $\nu_2(11664)=5$ and gcd$(5,18)=1$,  by Corollary \ref{nspe}, $K$ is monogenic and $\th$ generates a power integral basis.
    \item
    For $F(x)=x^{36}-37\in \Z[x]$, as $37\equiv 1\md{9}$, by Theorem \ref{npib}, $\Z_K$ is not monogenic.
    \item
    For $F(x)=x^{12}-13\in \Z[x]$,, since  $\nu_2=\nu_2(13-1)=2$, neither Theorem \ref{pib} nor Theorem \ref{npib}, can not give an answer about the monogenity of $K$.  Let us show that $2$ is a common index divisor of $K$ and so we conclude that $K$ is not monogenic. First $\ol{F(x)}=(x-1)^4(x^2+x+1)^4$ in $\F_2[x]$. For $\ph=x^2+x+1$,  
    $F(x)=\ph(x)^6+ (9-6x)\ph(x)^5-(25+5x)\ph(x)^4+(18+24x)\ph(x)^3 -18x\ph(x)^2 -(4-4x)\ph(x) -12$ is the $\ph$-expansion of $F(x)$. Thus with respect to $p=2$, $\npp{F}=S$ has  a single side joining the points $(0,2)$ and $(4,0)$ such that $F_S(y)=(1+x)y^2+xy+1$ in $\fph[y]$. First by Theorem \ref{ore},  $2$ divides the index $(\Z_K:\Z[\al])$.  But why  $2$ divides the index $(\Z_K:\Z[\th])$ for every  generator $\th\in \Z_K$ of $K$?  Since $F_S(y)=(y+1)(y+x^2)$ in $\fph[y]$, by Remark \ref{remore}, there are at least $2$ prime ideals of $\Z_K$ lying above $2$ with residue degree $2$ each. As $x^2+x+1$ is the unique monic irreducible polynomial over  $\F_2$, by Lemma \ref{index}, $2$ divides the index $(\Z_K:\Z[\th])$ for every  generator $\th\in \Z_K$ of $K$ and thus $K$ is not monogenic.
    \item
    For $F(x)=x^{12}-17\in \Z[x]$, since  $17\equiv -1\md9$, neither Theorem \ref{pib} nor Theorem \ref{npib}, can not give an answer about the monogenity of $K$.  Let us show that $3$ is a common index divisor of $K$ and so we conclude that $K$ is not monogenic. First $\ol{F(x)}=(x^2+x-1)^4(x^2-x-1)^4$ in $\F_3[x]$. Let $\ph_1=x^2+x-1$ and $\ph_2=x^2-x-1$. Then  
    $F(x)=\ph_1(x)^6+ (21-6x)\ph_1(x)^5+(125-65x)\ph_1(x)^4+(338-256x)\ph_1(x)^3 +(468-474x)\ph_1(x)^2 +(324-420x)\ph_1(x) +(72-144x)$ and 
    $F(x)=\ph_1(x)^6+ (21+6x)\ph_1(x)^5+(125+65x)\ph_1(x)^4+(338+256x)\ph_1(x)^3 +(468+474x)\ph_1(x)^2 +(324+420x)\ph_1(x) +(72+144x)$ are the $\ph_i$-expansion of $F(x)$ for $i=1,2$. Thus with respect to $p=3$, $N_{\ph_i}^+({F})=S_{i1}+S_{i2}$ has  two  sides  joining the points $(0,2)$, $(1,1)$, and $(4,0)$ (see $FIGURE 6$). Thus $d(S_{ij})=1$ for every $i,j=1,2$, and so every $\ph_i$ provides $2$ prime ideals of $\Z_K$ lying above $3$  with residue degree $2$ each.  It follows that there are $4$ prime ideals of $\Z_K$ lying above $3$ with residue degree $2$ each. Since there is only $3$ monic irreducible polynomial of degree $2$ over  $\F_3$, namely, $x^2$, $\ph_1$, and $\ph_2$, by Lemma \ref{index}, $3$ is a common index divisor of $K$ and $K$ is not monogenic.
 \centering
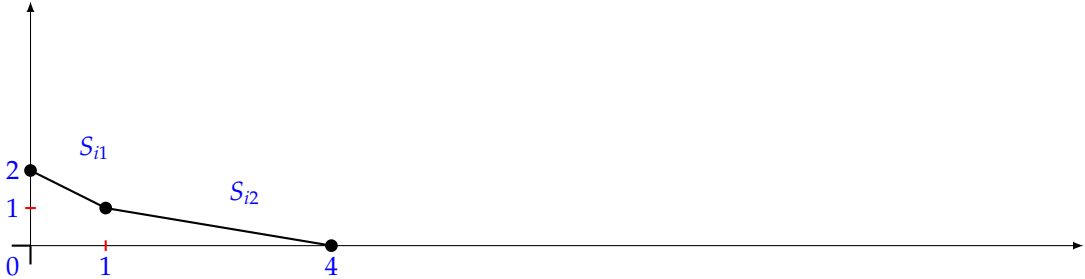
\begin{figure}[htbp] 
\begin{tikzpicture}[x=0.5cm,y=0.5cm]
\draw[latex-latex] (0,6.5) -- (0,0) -- (28,0) ;
\draw[thick] (0,0) -- (-0.5,0);
\draw[thick] (0,0) -- (0,-0.5); 
\draw[thick,red] (2,-2pt) -- (2,2pt);
\draw[thick,red] (8,-2pt) -- (8,2pt);
\draw[thick,red] (-2pt,1) -- (2pt,1);
\draw[thick,red] (-2pt,2) -- (2pt,2);
\node at (0,0) [below left,blue]{\footnotesize  $0$};
\node at (2,0) [below ,blue]{\footnotesize $1$};
\node at (8,0) [below ,blue]{\footnotesize  $4$};
\node at (0,1) [left ,blue]{\footnotesize  $1$};
\node at (0,2) [left ,blue]{\footnotesize  $2$};
\draw[thick,mark=*] plot coordinates{(0,2) (2,1) (8,0)};
\node at (1,2) [above right  ,blue]{\footnotesize  $S_{i1}$};
\node at (5,0.8) [above right  ,blue]{\footnotesize  $S_{i2}$};
\end{tikzpicture}
\caption{    \large   $N^+_{\ph_i}{F}$ \hspace{5cm}}
\end{figure}    
    \end{enumerate}
 	 

\begin{thebibliography}{99}
\bibitem{AN6} {\sc S. Ahmad, T. Nakahara, and S. M. Husnine}, {\em Power integral bases for certain pure sextic
fields}, Int. J. of Number Theory v:10, No 8 (2014) 2257--2265. 
\bibitem{AN} {\sc S. Ahmad, T. Nakahara, and A. Hameed}, {\em On certain pure sextic fields related to a problem of Hasse}, Int. J. Alg. and Comput. 26(3) (2016) 577--583 
\bibitem{Epr}
{\sc H. Ben Yakou  and L. El Fadil}, {\em on power integral bases  for certain pure  number fields defined by $x^{p^r}-m$},  Int. J. Number theory, 2021, doi:10.1142/S1793042121500858 
\bibitem{E2r5s} {\sc H. Ben yakkou, L. El Fadil, and A. Chillali}, {\em On Power integral bases for  certain  pure number fields defined by $x^{2^r\cdot 5^s}-m$}, Comm. in Algebra, 2021, doi: 10.1080/0092782.2021.1883642
		\bibitem{BK}
		{\sc H. Ben Yakou  and O. Kchit}, {\em on power integral bases  for certain pure  number fields defined by $x^{3^r}-m$}. São Paulo J. Math. Sci. (2021). https://doi.org/10.1007/s40863-021-00251-2

\bibitem{Co} H. Cohen, {\em A Course in Computational Algebraic Number Theory}, GTM 138, Springer-Verlag Berlin  Heidelberg (1993)
\bibitem{R} R. Dedekind,  {\em \"Uber den Zusammenhang zwischen der Theorie der Ideale und der Theorie der h\"oheren Kongruenzen}, G\"ottingen Abhandlungen {\bf 23} (1878) 1--23
\bibitem {E07} {\sc L. El Fadil}, {\em Computation of a power integral basis of a pure cubic number field}, Int. J. Contemp. Math. Sci. 2(13-16)(2007) 601--606
\bibitem{E6} {\sc L. El Fadil}, {\em On Power integral bases for certain pure sextic fields} (To appear in a forthcoming issue of Bol. Soc. Paran. Math.)
\bibitem{E6s} {\sc L. El Fadil}, {\em On Power integral bases for certain pure sextic fields}J, J. Number Theory, 228 (2021) 375–389
\bibitem{E12} L. El Fadil, {\em On Power integral bases for certain pure number fields } (To appear in a forthcoming issue of Publicationes Mathematicae)
\bibitem{E24}{\sc L. El Fadil}, {\em On Power integral bases for certain pure number fields defined by $x^{24}-m$ }, Stud. Sci. Math. Hung. {\bf 57}(3) (2020) 397–407
\bibitem{E18}{\sc L. El Fadil}, {\em On Power integral bases for certain pure number fields defined by $x^{18}-m$ }, (To appear in Comm. Math. Uni. of Carolina)
\bibitem{E23k}{\sc L. El Fadil}, 
{\em On Power integral bases for  certain  pure number fields defined by $x^{2\cdot 3^k}-m$} ( To appear in a forthcoming issue of Act. Arith.)
\bibitem{E3u7v} L. El Fadil, {\em On Power integral bases for certain pure number fields } (To appear in a forthcoming issue of Colloq. Math.)
\bibitem{El} L. El Fadil,{\em On Newton polygon's techniques and factorization of polynomial over Henselian valued fields},  J. of Algebra and  its Appl.  (2020), doi: S0219498820501881
\bibitem{EMN}{\sc L. El Fadil}, {\sc J. Montes} and 
{\sc E. Nart}, {\em Newton polygons and $p$-integral bases of quartic number fields}, 
J. Algebra and Appl. 11(4) (2012) 1250073
\bibitem{F4} {\sc T. Funakura}, {\em On integral bases of pure quartic fields}, Math. J.
Okayama Univ. 26 (1984) 27--41
\bibitem{G} {\sc  I. Ga\'al}, {\em Power integral bases in algebraic number fields}, Ann. Univ. Sci. Budapest. Sect. Comp. 18 (1999) 61--87 
\bibitem{G19} {\sc I. Ga\'al}, Diophantine equations and power integral bases, Theory and algorithm, Second edition, Boston, Birkh\"auser, 2019
\bibitem{GO} {\sc  I. Ga\'al, P. Olajos, and M. Pohst}, {\em Power integral bases in orders of composite fields}, Exp. Math. 11(1) (2002) 87–90. 
\bibitem{GR4} {\sc  I. Ga\'al and L. Remete}, {\em Binomial Thue equations and power integral bases
in  pure  quartic  fields}, JP  Journal  of  Algebra  Number  Theory  Appl. 32(1)
(2014) 49--61
\bibitem{GR17} {\sc  I. Ga\'al and L. Remete}, {\em Power integral bases and monogenity of pure  fields}, J. of Number Theory 173 (2017) 129--146 
\bibitem{GMN} J. Guardia, J. Montes and E. Nart,{\em Newton polygons of higher order in algebraic number
  theory}, J. trans. of ams  {\bf 364}(1) (2012)  361--416
  \bibitem{HN8} {\sc A.Hameed  and  T.Nakahara}, {\em Integral  bases  and  relative  monogenity  of pure octic fields}, Bull. Math. Soc. Sci. Math. R épub. Soc. Roum. 58(106) No. 4(2015) 419--433
\bibitem{CNS}  A. Hameed, T. Nakahara, S. M. Husnine, On existence of canonical number system in certain classes of pure algebraic number fields, J. Prime Res. Math. 7(2011) 19--24. 
{\bibitem{He2} {\sc K. Hensel}, {\em Arithmetische Untersuchungen \"{u}ber die gemeinsamen
				ausserwesentlichen Discriminantentheiler einer Gattung}, J. Reine Angew. Math., 113:128–160, 1894. ISSN 0075-4102.  doi: 10.1515/crll.1894.113.128.}
				\bibitem{MN}  J. Montes and E. Nart, {\em On a theorem of Ore},  J.  Algebra {\bf 146}(2) (1992) 318--334
 \bibitem{Neu}
{\sc J.  Neukirch}, {\em Algebraic Number Theory},
Springer-Verlag, Berlin, 1999.
\bibitem{O}{\sc O. Ore}, {\em Newtonsche
Polygone in der Theorie der algebraischen Korper}, Math. Ann., 99
(1928), 84--117
\bibitem{P}{\sc A. Peth\"o and M. Pohst},  {\em On the indices of multiquadratic number fields}, Acta Arith. 153(4) (2012) 393--414
	{\bibitem{HS}{\sc H. Smith}, {\em The monogeneity of radical extensions}. Acta Arith. 2021, 198(3): 313-327}.
\end{thebibliography}
\end{document}